\documentclass[a4paper,twoside]{article}
\usepackage{amsmath,amssymb,amsfonts}
\usepackage{hyperref}

\newtheorem{theorem}{Theorem}[section]
\newtheorem{proposition}[theorem]{Proposition}
\newtheorem{corollary}[theorem]{Corollary}
\newtheorem{lemma}[theorem]{Lemma}
\newtheorem{remark}[theorem]{Remark}

\input{tcilatex}
\begin{document}

\date{}
\title{\textbf{Bipolar-Valued Fuzzy Ideals in LA-semigroups}}
\author{\textbf{Naveed Yaqoob} \\
Department of Mathematics, Quaid-i-Azam University, Islamabad, Pakistan\\
nayaqoob@ymail.com\\
}
\maketitle

\begin{abstract}
In this paper, the notion of bipolar-valued fuzzy LA-subsemigroups is
introduced and also some properties of bipolar-valued fuzzy left (right,
bi-, interior) ideals of LA-semigroups has been discussed.
\end{abstract}

\emph{2010 AMS Classification}$:$ 06F35, 03E72, 20M12


\section{Introduction}

The concept of a fuzzy set was introduced by Zadeh \cite{zadeh}, in 1965.
Since its inception, the theory has developed in many directions and found
applications in a wide variety of fields. There has been a rapid growth in
the interest of fuzzy set theory and its applications from the past several
years. Many researchers published high-quality research articles on fuzzy
sets in a variety of international journals. The study of fuzzy set in
algebraic structure has been started in the definitive paper of Rosenfeld
1971 \cite{Rosen}. Fuzzy subgroup and its important properties were defined
and established by Rosenfeld \cite{Rosen}. In 1981, Kuroki introduced the
concept of fuzzy ideals and fuzzy bi-ideals in semigroups in his paper \cite%
{Kuroki}.

There are several kinds of fuzzy set extensions in the fuzzy set theory, for
example, intuitionistic fuzzy sets, interval-valued fuzzy sets, vague sets,
etc. Bipolar-valued fuzzy set is another extension of fuzzy set whose
membership degree range is different from the above extensions. Lee \cite%
{Lee} introduced the notion of bipolar-valued fuzzy sets. Bipolar-valued
fuzzy sets are an extension of fuzzy sets whose membership degree range is
enlarged from the interval $[0,1]$ to $[-1,1]$. In a bipolar-valued fuzzy
set, the membership degree $0$ indicate that elements are irrelevant to the
corresponding property, the membership degrees on $(0,1]$ assign that
elements somewhat satisfy the property, and the membership degrees on $%
[-1,0) $ assign that elements somewhat satisfy the implicit counter-property 
\cite{Lee, Lee2}.

Akram et al. \cite{Akram} introduced the concept of bipolar fuzzy
K-algebras. In \cite{Jun}, Jun and park applied the notion of bipolar-valued
fuzzy sets to BCH-algebras. They introduced the concept of bipolar fuzzy
subalgebras and bipolar fuzzy ideals of a BCH-algebra. Lee \cite{KLee2}
applied the notion of bipolar fuzzy subalgebras and bipolar fuzzy ideals of
BCK/BCI-algebras. Also some results on bipolar-valued fuzzy BCK/BCI-algebras
are introduced by Saeid in \cite{Saeid}.

This paper concerns the relationship between bipolar-valued fuzzy sets and
left almost semigroups. The left almost semigroup abbreviated as an
LA-semigroup, was first introduced by Kazim and Naseerudin \cite{Kazim}.
They generalized some useful results of semigroup theory. They introduced
braces on the left of the ternary commutative law $abc=cba,$ to get a new
pseudo associative law, that is $(ab)c=(cb)a,$ and named it as left
invertive law. An LA-semigroup is the midway structure between a commutative
semigroup and a groupoid. Despite the fact, the structure is non-associative
and non-commutative. It nevertheless possesses many interesting properties
which we usually find in commutative and associative algebraic structures.
Mushtaq and Yusuf produced useful results \cite{1979}, on locally
associative LA-semigroups in 1979. In this structure they defined powers of
an element and congruences using these powers. They constructed quotient
LA-semigroups using these congruences.\ It is a useful non-associative
structure with wide applications in theory of flocks.

In this paper, we have introduced the notion of bipolar-valued fuzzy
LA-subsemigroups and bipolar-valued fuzzy left (right, bi-, interior) ideals
in LA-semigroups.

\section{\textbf{Preliminaries and basic definitions}}

\textbf{Definition 2.1. }\cite{Kazim} A groupoid $\left( S,\cdot \right) $\
is called an LA-semigroup$,$ if it satisfies left invertive law 
\begin{equation*}
\left( a\cdot b\right) \cdot c=\left( c\cdot b\right) \cdot a,\text{ \ for
all }a,b,c\in S\text{.}
\end{equation*}

\textbf{Example 2.1 }\cite{1978} Let $\left( 
\mathbb{Z}
,+\right) $\ denote the commutative group of integers under addition. Define
a binary operation \textquotedblleft $\ast $\textquotedblright\ in $%
\mathbb{Z}
$\ as follows$:$%
\begin{equation*}
a\ast b=b-a,\text{ \ for all }a,b\in 
\mathbb{Z}
\text{.}
\end{equation*}%
Where \textquotedblleft $-$\textquotedblright\ denotes the ordinary
subtraction of integers. Then $\left( 
\mathbb{Z}
,\ast \right) $\ is an LA-semigroup.

\textbf{Example 2.2 }\cite{1978} Define a binary operation \textquotedblleft 
$\ast $\textquotedblright\ in $%
\mathbb{R}
$\ as follows$:$%
\begin{equation*}
a\ast b=b\div a,\text{ \ for all }a,b\in 
\mathbb{R}
\text{.}
\end{equation*}%
Then $\left( 
\mathbb{R}
,\ast \right) $\ is an LA-semigroup.

\begin{lemma}
\label{l1}\cite{1979} If $S$ is an LA-semigroup with left identity $e$, then 
$a(bc)=b(ac)$ for all $a,b,c\in S.$
\end{lemma}

Let $S$ be an LA-semigroup. A nonempty subset $A$ of $S$ is called an
LA-subsemigroup of $S$ if $ab\in A$ for all $a,b\in A$. A nonempty subset $L$
of $S$ is called a left ideal of $S$ if $SL\subseteq L$ and a nonempty
subset $R$ of $S$ is called a right ideal of $S$ if $RS\subseteq R$. A
nonempty subset $I$ of $S$ is called an ideal of $S$ if $I$ is both a left
and a right ideal of $S$. A subset $A$ of $S$ is called an interior ideal of 
$S$ if $(SA)S\subseteq A$. An LA-subsemigroup $A$ of $S$ is called a
bi-ideal of $S$ if $(AS)A\subseteq A$.

In an LA-semigroup the medial law holds:%
\begin{equation*}
(ab)(cd)=(ac)(bd),\text{ \ for all }a,b,c,d\in S.
\end{equation*}

In an LA-semigroup $S$ with left identity, the paramedial law holds:%
\begin{equation*}
(ab)(cd)=(dc)(ba),\text{ \ for all }a,b,c,d\in S.
\end{equation*}

Now we will recall the concept of bipolar-valued fuzzy sets.

\textbf{Definition 2.2 }\cite{Lee2} Let $X$\ be a nonempty set. A
bipolar-valued fuzzy subset (BVF-subset, in short) $B$\ of $X$\ is an object
having the form%
\begin{equation*}
B=\left\{ \left\langle x,\mu _{B}^{+}(x),\mu _{B}^{-}(x)\right\rangle :x\in
X\right\} .
\end{equation*}%
Where $\mu _{B}^{+}:X\rightarrow \lbrack 0,1]$\ and $\mu
_{B}^{-}:X\rightarrow \lbrack -1,0]$.

The positive membership degree $\mu _{B}^{+}(x)$ denotes the satisfaction
degree of an element $x$ to the property corresponding to a bipolar-valued
fuzzy set $B=\left\{ \left\langle x,\mu _{B}^{+}(x),\mu
_{B}^{-}(x)\right\rangle :x\in X\right\} $, and the negative membership
degree $\mu _{B}^{-}(x)$ denotes the satisfaction degree of $x$ to some
implicit counter property of $B=\left\{ \left\langle x,\mu _{B}^{+}(x),\mu
_{B}^{-}(x)\right\rangle :x\in X\right\} $. For the sake of simplicity, we
shall use the symbol $B=\left\langle \mu _{B}^{+},\mu _{B}^{-}\right\rangle $
for the bipolar-valued fuzzy set $B=\left\{ \left\langle x,\mu
_{B}^{+}(x),\mu _{B}^{-}(x)\right\rangle :x\in X\right\} .$

\textbf{Definition 2.3 }Let $B_{1}=\left\langle \mu _{B_{1}}^{+},\mu
_{B_{1}}^{-}\right\rangle $ and $B_{2}=\left\langle \mu _{B_{2}}^{+},\mu
_{B_{2}}^{-}\right\rangle $ be two BVF-subsets of a nonempty set $X$. Then
the product of two BVF-subsets is denoted by $B_{1}\circ B_{2}$ and defined
as: 
\begin{eqnarray*}
\left( \mu _{B_{1}}^{+}\circ \mu _{B_{2}}^{+}\right) \left( x\right) 
&=&\left\{ 
\begin{array}{l}
\tbigvee_{x=yz}\left\{ \mu _{B_{1}}^{+}\left( y\right) \wedge \mu
_{B_{2}}^{+}\left( z\right) \right\} ,\text{ if }x=yz\text{ for some }y,z\in
S \\ 
\text{ }0\text{ \ \ \ \ \ \ \ \ \ \ \ \ \ \ \ \ \ \ \ \ \ \ \ \ \ \ \ \ \ \
\ otherwise.\ \ \ \ \ \ \ \ \ \ \ \ \ \ \ \ \ \ \ \ \ \ \ \ \ \ \ \ \ \ \ }%
\end{array}%
\right.  \\
\left( \mu _{B_{1}}^{-}\circ \mu _{B_{2}}^{-}\right) \left( x\right) 
&=&\left\{ 
\begin{array}{l}
\tbigwedge_{x=yz}\left\{ \mu _{B_{1}}^{-}\left( y\right) \vee \mu
_{B_{2}}^{-}\left( z\right) \right\} ,\text{ if }x=yz\text{ for some }y,z\in
S \\ 
\text{ }0\text{ \ \ \ \ \ \ \ \ \ \ \ \ \ \ \ \ \ \ \ \ \ \ \ \ \ \ \ \ \ \
\ otherwise.\ \ \ \ \ \ \ \ \ \ \ \ \ \ \ \ \ \ \ \ \ \ \ \ \ \ \ \ \ \ \ }%
\end{array}%
\right. 
\end{eqnarray*}

Note that an LA-semigroup $S$ can be considered as a BVF-subset of itself
and let%
\begin{eqnarray*}
\Gamma &=&\left\langle \mathcal{S}_{\Gamma }^{+}(x),\mathcal{S}_{\Gamma
}^{-}(x)\right\rangle \\
&=&\left\{ \left\langle x,\mathcal{S}_{\Gamma }^{+}(x),\mathcal{S}_{\Gamma
}^{-}(x)\right\rangle :\mathcal{S}_{\Gamma }^{+}(x)=1\text{ and }\mathcal{S}%
_{\Gamma }^{-}(x)=-1,\text{ for all }x\text{ in }S\right\}
\end{eqnarray*}%
be a BVF-subset and $\Gamma =\left\langle \mathcal{S}_{\Gamma }^{+}(x),%
\mathcal{S}_{\Gamma }^{-}(x)\right\rangle $ will be carried out in
operations with a BVF-subset $B=\left\langle \mu _{B}^{+},\mu
_{B}^{-}\right\rangle $ such that $\mathcal{S}_{\Gamma }^{+}$ and $\mathcal{S%
}_{\Gamma }^{-}$ will be used in collaboration with $\mu _{B}^{+}$ and $\mu
_{B}^{-}$ respectively.

Let $BVF(S)$ denote the set of all BVF-subsets of an LA-semigroup $S.$

\begin{proposition}
\label{P1}Let $S$ be an LA-semigroup, then the set $(BVF(S),\circ )$ is an
LA-semigroup.
\end{proposition}

\textbf{Proof. }Clearly $BVF(S)$ is closed. Let $B_{1}=\left\langle \mu
_{B_{1}}^{+},\mu _{B_{1}}^{-}\right\rangle ,$ $B_{2}=\left\langle \mu
_{B_{2}}^{+},\mu _{B_{2}}^{-}\right\rangle $ and $B_{3}=\left\langle \mu
_{B_{3}}^{+},\mu _{B_{3}}^{-}\right\rangle $ be in $BVF(S)$. Let $x$ be any
element of $S$ such that $x\neq yz$ for some $y,z\in S$. Then we have%
\begin{equation*}
\left( \left( \mu _{B_{1}}^{+}\circ \mu _{B_{2}}^{+}\right) \circ \mu
_{B_{3}}^{+}\right) (x)=0=\left( \left( \mu _{B_{3}}^{+}\circ \mu
_{B_{2}}^{+}\right) \circ \mu _{B_{1}}^{+}\right) (x).
\end{equation*}%
And%
\begin{equation*}
\left( \left( \mu _{B_{1}}^{-}\circ \mu _{B_{2}}^{-}\right) \circ \mu
_{B_{3}}^{-}\right) (x)=0=\left( \left( \mu _{B_{3}}^{-}\circ \mu
_{B_{2}}^{-}\right) \circ \mu _{B_{1}}^{-}\right) (x).
\end{equation*}%
Let $x$ be any element of $S$ such that $x=yz$ for some $y,z\in S$. Then we
have%
\begin{eqnarray*}
\left( \left( \mu _{B_{1}}^{+}\circ \mu _{B_{2}}^{+}\right) \circ \mu
_{B_{3}}^{+}\right) (x) &=&{\bigvee }_{x=yz}\left\{ \left( \mu
_{B_{1}}^{+}\circ \mu _{B_{2}}^{+}\right) (y)\wedge \mu
_{B_{3}}^{+}(z)\right\} \\
&=&{\bigvee }_{x=yz}\left\{ \left( {\bigvee }_{y=pq}\left\{ \mu
_{B_{1}}^{+}(p)\wedge \mu _{B_{2}}^{+}(q)\right\} \right) \wedge \mu
_{B_{3}}^{+}(z)\right\} \\
&=&{\bigvee }_{x=yz}\text{ }{\bigvee }_{y=pq}\left\{ \mu
_{B_{1}}^{+}(p)\wedge \mu _{B_{2}}^{+}(q)\wedge \mu _{B_{3}}^{+}(z)\right\}
\\
&=&{\bigvee }_{x=(pq)z}\left\{ \mu _{B_{1}}^{+}(p)\wedge \mu
_{B_{2}}^{+}(q)\wedge \mu _{B_{3}}^{+}(z)\right\} \\
&=&{\bigvee }_{x=(zq)p}\left\{ \mu _{B_{3}}^{+}(z)\wedge \mu
_{B_{2}}^{+}(q)\wedge \mu _{B_{1}}^{+}(p)\right\} \\
&=&{\bigvee }_{x=sp}\left\{ \left( {\bigvee }_{s=zq}\left\{ \mu
_{B_{3}}^{+}(z)\wedge \mu _{B_{2}}^{+}(q)\right\} \right) \wedge \mu
_{B_{1}}^{+}(p)\right\} \\
&=&{\bigvee }_{x=sp}\left\{ \left( \mu _{B_{3}}^{+}\circ \mu
_{B_{2}}^{+}\right) (s)\wedge \mu _{B_{1}}^{+}(p)\right\} \\
&=&\left( \left( \mu _{B_{3}}^{+}\circ \mu _{B_{2}}^{+}\right) \circ \mu
_{B_{1}}^{+}\right) (x).
\end{eqnarray*}%
And%
\begin{eqnarray*}
\left( \left( \mu _{B_{1}}^{-}\circ \mu _{B_{2}}^{-}\right) \circ \mu
_{B_{3}}^{-}\right) (x) &=&{\bigwedge }_{x=yz}\left\{ \left( \mu
_{B_{1}}^{-}\circ \mu _{B_{2}}^{-}\right) (y)\vee \mu _{B_{3}}^{-}(z)\right\}
\\
&=&{\bigwedge }_{x=yz}\left\{ \left( {\bigwedge }_{y=pq}\left\{ \mu
_{B_{1}}^{-}(p)\vee \mu _{B_{2}}^{-}(q)\right\} \right) \vee \mu
_{B_{3}}^{-}(z)\right\} \\
&=&{\bigwedge }_{x=yz}\text{ }{\bigwedge }_{y=pq}\left\{ \mu
_{B_{1}}^{-}(p)\vee \mu _{B_{2}}^{-}(q)\vee \mu _{B_{3}}^{-}(z)\right\} \\
&=&{\bigwedge }_{x=(pq)z}\left\{ \mu _{B_{1}}^{-}(p)\vee \mu
_{B_{2}}^{-}(q)\vee \mu _{B_{3}}^{-}(z)\right\} \\
&=&{\bigwedge }_{x=(zq)p}\left\{ \mu _{B_{3}}^{-}(z)\vee \mu
_{B_{2}}^{-}(q)\vee \mu _{B_{1}}^{-}(p)\right\} \\
&=&{\bigwedge }_{x=sp}\left\{ \left( {\bigwedge }_{s=zq}\left\{ \mu
_{B_{3}}^{-}(z)\vee \mu _{B_{2}}^{-}(q)\right\} \right) \vee \mu
_{B_{1}}^{-}(p)\right\} \\
&=&{\bigwedge }_{x=sp}\left\{ \left( \mu _{B_{3}}^{-}\circ \mu
_{B_{2}}^{-}\right) (s)\vee \mu _{B_{1}}^{-}(p)\right\} \\
&=&\left( \left( \mu _{B_{3}}^{-}\circ \mu _{B_{2}}^{-}\right) \circ \mu
_{B_{1}}^{-}\right) (x).
\end{eqnarray*}

Hence $(BVF(S),\circ )$ is an LA-semigroup. $\ \ \Box $

\begin{corollary}
\label{C1}If $S$ is an LA-semigroup, then the medial law holds in $BVF(S)$.
\end{corollary}

\textbf{Proof. }Let $B_{1}=\left\langle \mu _{B_{1}}^{+},\mu
_{B_{1}}^{-}\right\rangle $, $B_{2}=\left\langle \mu _{B_{2}}^{+},\mu
_{B_{2}}^{-}\right\rangle ,$ $B_{3}=\left\langle \mu _{B_{3}}^{+},\mu
_{B_{3}}^{-}\right\rangle $ and $B_{4}=\left\langle \mu _{B_{4}}^{+},\mu
_{B_{4}}^{-}\right\rangle $ be in $BVF(S)$. By successive use of left
invertive law%
\begin{eqnarray*}
\left( \mu _{B_{1}}^{+}\circ \mu _{B_{2}}^{+}\right) \circ \left( \mu
_{B_{3}}^{+}\circ \mu _{B_{4}}^{+}\right) &=&\left( \left( \mu
_{B_{3}}^{+}\circ \mu _{B_{4}}^{+}\right) \circ \mu _{B_{2}}^{+}\right)
\circ \mu _{B_{1}}^{+} \\
&=&\left( \left( \mu _{B_{2}}^{+}\circ \mu _{B_{4}}^{+}\right) \circ \mu
_{B_{3}}^{+}\right) \circ \mu _{B_{1}}^{+} \\
&=&\left( \mu _{B_{1}}^{+}\circ \mu _{B_{3}}^{+}\right) \circ \left( \mu
_{B_{2}}^{+}\circ \mu _{B_{4}}^{+}\right) .
\end{eqnarray*}%
And%
\begin{eqnarray*}
\left( \mu _{B_{1}}^{-}\circ \mu _{B_{2}}^{-}\right) \circ \left( \mu
_{B_{3}}^{-}\circ \mu _{B_{4}}^{-}\right) &=&\left( \left( \mu
_{B_{3}}^{-}\circ \mu _{B_{4}}^{-}\right) \circ \mu _{B_{2}}^{-}\right)
\circ \mu _{B_{1}}^{-} \\
&=&\left( \left( \mu _{B_{2}}^{-}\circ \mu _{B_{4}}^{-}\right) \circ \mu
_{B_{3}}^{-}\right) \circ \mu _{B_{1}}^{-} \\
&=&\left( \mu _{B_{1}}^{-}\circ \mu _{B_{3}}^{-}\right) \circ \left( \mu
_{B_{2}}^{-}\circ \mu _{B_{4}}^{-}\right) .
\end{eqnarray*}%
Hence this shows that the medial law holds in $BVF(S)$. $\ \ \Box $

\section{\textbf{Bipolar-valued fuzzy ideals in LA-semigroup}}

\textbf{Definition 3.1 }A BVF-subset $B=\left\langle \mu _{B}^{+},\mu
_{B}^{-}\right\rangle $ of an LA-semigroup $S$ is called a bipolar-valued
fuzzy LA-subsemigroup of $S$ if%
\begin{equation*}
\mu _{B}^{+}\left( xy\right) \geq \mu _{B}^{+}\left( x\right) \wedge \mu
_{B}^{+}\left( y\right) \text{\ \ and \ \ }\mu _{B}^{-}\left( xy\right) \leq
\mu _{B}^{-}\left( x\right) \vee \mu _{B}^{-}\left( y\right)
\end{equation*}%
for all $x,y\in S$.

\textbf{Definition 3.2 }A BVF-subset $B=\left\langle \mu _{B}^{+},\mu
_{B}^{-}\right\rangle $ of an LA-semigroup $S$ is called a bipolar-valued
fuzzy left ideal of $S$ if%
\begin{equation*}
\mu _{B}^{+}\left( xy\right) \geq \mu _{B}^{+}\left( y\right) \text{ \ \ and
\ \ }\mu _{B}^{-}\left( xy\right) \leq \mu _{B}^{-}\left( y\right)
\end{equation*}%
for all $x,y\in S$.

\textbf{Definition 3.3 }A BVF-subset $B=\left\langle \mu _{B}^{+},\mu
_{B}^{-}\right\rangle $ of an LA-semigroup $S$ is called a bipolar-valued
fuzzy right ideal of $S$ if%
\begin{equation*}
\mu _{B}^{+}\left( xy\right) \geq \mu _{B}^{+}\left( x\right) \text{ \ \ and
\ }\mu _{B}^{-}\left( xy\right) \leq \mu _{B}^{-}\left( x\right)
\end{equation*}%
for all $x,y\in S$.

A BVF-subset $B=\left\langle \mu _{B}^{+},\mu _{B}^{-}\right\rangle $ of an
LA-semigroup $S$ is called a BVF-ideal or BVF-two-sided ideal of $S$ if $%
B=\left\langle \mu _{B}^{+},\mu _{B}^{-}\right\rangle $ is both BVF-left and
BVF-right ideal of $S$.

\textbf{Example 3.1} Let $S=\{a,b,c,d\}$, the binary operation "$\cdot $" on 
$S$ be defined as follows:%
\begin{equation*}
\begin{tabular}{l|llll}
$\cdot $ & $a$ & $b$ & $c$ & $d$ \\ \hline
$a$ & $b$ & $d$ & $c$ & $a$ \\ 
$b$ & $a$ & $b$ & $c$ & $d$ \\ 
$c$ & $c$ & $c$ & $c$ & $c$ \\ 
$d$ & $d$ & $a$ & $c$ & $b$%
\end{tabular}%
\end{equation*}%
Clearly, $S$ is an LA-semigroup. But $S$ is not a semigroup because $%
d=d\cdot (b\cdot a)\neq (d\cdot b)\cdot a=b.$ Now we define BVF-subset as%
\begin{equation*}
B=\left\langle \mu _{B}^{+},\mu _{B}^{-}\right\rangle =\left\langle \left( 
\frac{a}{0.2},\frac{b}{0.2},\frac{c}{0.7},\frac{d}{0.2}\right) ,\text{ }%
\left( \frac{a}{-0.5},\frac{b}{-0.5},\frac{c}{-0.8},\frac{d}{-0.5}\right)
\right\rangle .
\end{equation*}%
Clearly $B$ is a BVF-ideal of $S.$

\begin{proposition}
\label{BP2}Every BVF-left (BVF-right) ideal $B=\left\langle \mu _{B}^{+},\mu
_{B}^{-}\right\rangle $ of an LA-semigroup $S$ is a bipolar-valued fuzzy
LA-subsemigroup of $S$.
\end{proposition}

\textbf{Proof. }Let $B=\left\langle \mu _{B}^{+},\mu _{B}^{-}\right\rangle $
be a BVF-left ideal of $S$ and for any $x,y\in S$,

\begin{equation*}
\mu _{B}^{+}\left( xy\right) \geq \mu _{B}^{+}\left( y\right) \geq \mu
_{B}^{+}\left( x\right) \wedge \mu _{B}^{+}\left( y\right) \text{.}
\end{equation*}%
And%
\begin{equation*}
\mu _{B}^{-}\left( xy\right) \leq \mu _{B}^{-}\left( y\right) \leq \mu
_{B}^{-}\left( x\right) \vee \mu _{B}^{-}\left( y\right) .
\end{equation*}%
Hence $B=\left\langle \mu _{B}^{+},\mu _{B}^{-}\right\rangle $ is a
bipolar-valued fuzzy LA-subsemigroup of $S$. The other case can be prove in
a similar way. $\ \ \Box $

\begin{lemma}
\label{L56}Let $B=\left\langle \mu _{B}^{+},\mu _{B}^{-}\right\rangle $ be a
BVF-subset of an LA-semigroup $S.$ Then
\end{lemma}

(1) $B=\left\langle \mu _{B}^{+},\mu _{B}^{-}\right\rangle $ is a
BVF-LA-subsemigroup of $S$\ if and only if $\mu _{B}^{+}\circ \mu
_{B}^{+}\subseteq \mu _{B}^{+}$ and $\mu _{B}^{-}\circ \mu _{B}^{-}\supseteq
\mu _{B}^{-}.$

(2) $B=\left\langle \mu _{B}^{+},\mu _{B}^{-}\right\rangle $ is a BVF-left
(resp. BVF-right) ideal of $S$\ if and only if $\mathcal{S}_{\Gamma
}^{+}\circ \mu _{B}^{+}\subseteq \mu _{B}^{+}$ and $\mathcal{S}_{\Gamma
}^{-}\circ \mu _{B}^{-}\supseteq \mu _{B}^{-}$ (resp. $\mu _{B}^{+}\circ 
\mathcal{S}_{\Gamma }^{+}\subseteq \mu _{B}^{+}$ and $\mu _{B}^{-}\circ 
\mathcal{S}_{\Gamma }^{-}\supseteq \mu _{B}^{-}$)$.$

\textbf{Proof. }(1)\textbf{\ }Let $B=\left\langle \mu _{B}^{+},\mu
_{B}^{-}\right\rangle $ be a BVF-LA-subsemigroup of $S$ and $x\in S.$ If $%
\left( \mu _{B}^{+}\circ \mu _{B}^{+}\right) \left( x\right) =0$ and $\left(
\mu _{B}^{-}\circ \mu _{B}^{-}\right) \left( x\right) =0,$ then $\left( \mu
_{B}^{+}\circ \mu _{B}^{+}\right) \left( x\right) \leq \mu _{B}^{+}\left(
x\right) $ and $\left( \mu _{B}^{-}\circ \mu _{B}^{-}\right) \left( x\right)
\geq \mu _{B}^{-}\left( x\right) .$ Otherwise,%
\begin{equation*}
\left( \mu _{B}^{+}\circ \mu _{B}^{+}\right) \left( x\right) ={\bigvee }%
_{x=yz}\left\{ \mu _{B}^{+}\left( y\right) \wedge \mu _{B}^{+}\left(
z\right) \right\} \leq {\bigvee }_{x=yz}\mu _{B}^{+}\left( yz\right) =\mu
_{B}^{+}\left( x\right) .
\end{equation*}%
And%
\begin{equation*}
\left( \mu _{B}^{-}\circ \mu _{B}^{-}\right) \left( x\right) ={\bigwedge }%
_{x=yz}\left\{ \mu _{B}^{-}\left( y\right) \vee \mu _{B}^{-}\left( z\right)
\right\} \geq {\bigwedge }_{x=yz}\mu _{B}^{-}\left( yz\right) =\mu
_{B}^{-}\left( x\right) .
\end{equation*}%
Thus $\mu _{B}^{+}\circ \mu _{B}^{+}\subseteq \mu _{B}^{+}$ and $\mu
_{B}^{-}\circ \mu _{B}^{-}\supseteq \mu _{B}^{-}.$\newline
Conversely, let $\mu _{B}^{+}\circ \mu _{B}^{+}\subseteq \mu _{B}^{+}$, $\mu
_{B}^{-}\circ \mu _{B}^{-}\supseteq \mu _{B}^{-}$ and $x,y\in S,$ then%
\begin{equation*}
\mu _{B}^{+}\left( xy\right) \geq \left( \mu _{B}^{+}\circ \mu
_{B}^{+}\right) \left( xy\right) ={\bigvee }_{xy=ab}\left\{ \mu
_{B}^{+}\left( a\right) \wedge \mu _{B}^{+}\left( b\right) \right\} \geq \mu
_{B}^{+}\left( x\right) \wedge \mu _{B}^{+}\left( y\right) .
\end{equation*}%
And%
\begin{equation*}
\mu _{B}^{-}\left( xy\right) \leq \left( \mu _{B}^{-}\circ \mu
_{B}^{-}\right) \left( xy\right) ={\bigwedge }_{xy=ab}\left\{ \mu
_{B}^{-}\left( a\right) \vee \mu _{B}^{-}\left( b\right) \right\} \leq \mu
_{B}^{-}\left( x\right) \vee \mu _{B}^{-}\left( y\right) .
\end{equation*}%
So $B=\left\langle \mu _{B}^{+},\mu _{B}^{-}\right\rangle $ is a
BVF-LA-subsemigroup of $S.$

(2) Let $B=\left\langle \mu _{B}^{+},\mu _{B}^{-}\right\rangle $ be a
BVF-left ideal of $S$ and $x\in S.$ If $\left( \mathcal{S}_{\Gamma
}^{+}\circ \mu _{B}^{+}\right) \left( x\right) =0$ and $\left( \mathcal{S}%
_{\Gamma }^{-}\circ \mu _{B}^{-}\right) \left( x\right) =0,$ then $\left( 
\mathcal{S}_{\Gamma }^{+}\circ \mu _{B}^{+}\right) \left( x\right) \leq \mu
_{B}^{+}\left( x\right) $ and $\left( \mathcal{S}_{\Gamma }^{-}\circ \mu
_{B}^{-}\right) \left( x\right) \geq \mu _{B}^{-}\left( x\right) .$
Otherwise,%
\begin{eqnarray*}
\left( \mathcal{S}_{\Gamma }^{+}\circ \mu _{B}^{+}\right) \left( x\right) &=&%
{\bigvee }_{x=ab}\left\{ \mathcal{S}_{\Gamma }^{+}\left( a\right) \wedge \mu
_{B}^{+}\left( b\right) \right\} ={\bigvee }_{x=ab}\left\{ 1\wedge \mu
_{B}^{+}\left( b\right) \right\} \\
&=&{\bigvee }_{x=ab}\mu _{B}^{+}\left( b\right) \leq {\bigvee }_{x=ab}\mu
_{B}^{+}\left( ab\right) =\mu _{B}^{+}\left( x\right) .
\end{eqnarray*}%
And%
\begin{eqnarray*}
\left( \mathcal{S}_{\Gamma }^{-}\circ \mu _{B}^{-}\right) \left( x\right) &=&%
{\bigwedge }_{x=ab}\left\{ \mathcal{S}_{\Gamma }^{-}\left( a\right) \vee \mu
_{B}^{-}\left( b\right) \right\} ={\bigwedge }_{x=ab}\left\{ -1\vee \mu
_{B}^{-}\left( b\right) \right\} \\
&=&{\bigwedge }_{x=ab}\mu _{B}^{-}\left( b\right) \geq {\bigwedge }%
_{x=ab}\mu _{B}^{-}\left( ab\right) =\mu _{B}^{-}\left( x\right) .
\end{eqnarray*}%
Thus $\mathcal{S}_{\Gamma }^{+}\circ \mu _{B}^{+}\subseteq \mu _{B}^{+}$ and 
$\mathcal{S}_{\Gamma }^{-}\circ \mu _{B}^{-}\supseteq \mu _{B}^{-}.$\newline
Conversely, let $\mathcal{S}_{\Gamma }^{+}\circ \mu _{B}^{+}\subseteq \mu
_{B}^{+}$, $\mathcal{S}_{\Gamma }^{-}\circ \mu _{B}^{-}\supseteq \mu
_{B}^{-} $ and $x,y\in S,$ then 
\begin{eqnarray*}
\mu _{B}^{+}\left( xy\right) &\geq &\left( \mathcal{S}_{\Gamma }^{+}\circ
\mu _{B}^{+}\right) \left( xy\right) ={\bigvee }_{xy=ab}\left\{ \mathcal{S}%
_{\Gamma }^{+}\left( a\right) \wedge \mu _{B}^{+}\left( b\right) \right\} \\
&\geq &\mathcal{S}_{\Gamma }^{+}\left( x\right) \wedge \mu _{B}^{+}\left(
y\right) =1\wedge \mu _{B}^{+}\left( y\right) =\mu _{B}^{+}\left( y\right) .
\end{eqnarray*}%
And 
\begin{eqnarray*}
\mu _{B}^{-}\left( xy\right) &\leq &\left( \mathcal{S}_{\Gamma }^{-}\circ
\mu _{B}^{-}\right) \left( xy\right) ={\bigwedge }_{xy=ab}\left\{ \mathcal{S}%
_{\Gamma }^{-}\left( a\right) \vee \mu _{B}^{-}\left( b\right) \right\} \\
&\leq &\mathcal{S}_{\Gamma }^{-}\left( x\right) \vee \mu _{B}^{-}\left(
y\right) =-1\vee \mu _{B}^{-}\left( y\right) =\mu _{B}^{-}\left( y\right) .
\end{eqnarray*}%
Thus $\mu _{B}^{+}\left( xy\right) \geq \mu _{B}^{+}\left( y\right) $ and $%
\mu _{B}^{-}\left( xy\right) \leq \mu _{B}^{-}\left( y\right) .$\ Thus $%
B=\left\langle \mu _{B}^{+},\mu _{B}^{-}\right\rangle $ is a BVF-left ideal
of $S.$ The second case can be seen in a similar way. $\ \ \Box $

Let $B_{1}=\left\langle \mu _{B_{1}}^{+},\mu _{B_{1}}^{-}\right\rangle $ and 
$B_{2}=\left\langle \mu _{B_{2}}^{+},\mu _{B_{2}}^{-}\right\rangle $ be two
BVF-subsets of an LA-semigroup $S.$ The symbol $B_{1}\cap B_{2}$ will mean
the following

\begin{equation*}
\left( \mu _{B_{1}}^{+}\cap \mu _{B_{2}}^{+}\right) (x)=\mu
_{B_{1}}^{+}(x)\wedge \mu _{B_{2}}^{+}(x),\text{ for all }x\in S.
\end{equation*}

\begin{equation*}
\left( \mu _{B_{1}}^{-}\cup \mu _{B_{2}}^{-}\right) (x)=\mu
_{B_{1}}^{-}(x)\vee \mu _{B_{2}}^{-}(x),\text{ for all }x\in S.
\end{equation*}%
The symbol $A\cup B$ will mean the following

\begin{equation*}
\left( \mu _{B_{1}}^{+}\cup \mu _{B_{2}}^{+}\right) (x)=\mu
_{B_{1}}^{+}(x)\vee \mu _{B_{2}}^{+}(x),\text{ for all }x\in S.
\end{equation*}

\begin{equation*}
\left( \mu _{B_{1}}^{-}\cap \mu _{B_{2}}^{-}\right) (x)=\mu
_{B_{1}}^{-}(x)\wedge \mu _{B_{2}}^{-}(x),\text{ for all }x\in S.
\end{equation*}

\begin{theorem}
Let $S$ be an LA-semigroup and $B_{1}=\left\langle \mu _{B_{1}}^{+},\mu
_{B_{1}}^{-}\right\rangle $ be a BVF-right ideal of $S$ and $%
B_{2}=\left\langle \mu _{B_{2}}^{+},\mu _{B_{2}}^{-}\right\rangle $ be a
BVF-left ideal of $S,$ then $B_{1}\circ B_{2}\subseteq B_{1}\cap B_{2}.$
\end{theorem}

\textbf{Proof. }Let for any $x,y,z\in S,$ if $x\neq yz,$ then we have%
\begin{equation*}
\left( \mu _{B_{1}}^{+}\circ \mu _{B_{2}}^{+}\right) (x)=0\leq \mu
_{B_{1}}^{+}(x)\wedge \mu _{B_{2}}^{+}(x)=\left( \mu _{B_{1}}^{+}\cap \mu
_{B_{2}}^{+}\right) (x).
\end{equation*}%
And%
\begin{equation*}
\left( \mu _{B_{1}}^{-}\circ \mu _{B_{2}}^{-}\right) (x)=0\geq \mu
_{B_{1}}^{-}(x)\vee \mu _{B_{2}}^{-}(x)=\left( \mu _{B_{1}}^{-}\cup \mu
_{B_{2}}^{-}\right) (x).
\end{equation*}%
Otherwise%
\begin{eqnarray*}
\left( \mu _{B_{1}}^{+}\circ \mu _{B_{2}}^{+}\right) (x) &=&{\bigvee }%
_{x=yz}\left\{ \mu _{B_{1}}^{+}\left( y\right) \wedge \mu _{B_{2}}^{+}\left(
z\right) \right\}  \\
&\leq &{\bigvee }_{x=yz}\left\{ \mu _{B_{1}}^{+}\left( yz\right) \wedge \mu
_{B_{2}}^{+}\left( yz\right) \right\}  \\
&=&\left\{ \mu _{B_{1}}^{+}\left( x\right) \wedge \mu _{B_{2}}^{+}\left(
x\right) \right\}  \\
&=&\left( \mu _{B_{1}}^{+}\cap \mu _{B_{2}}^{+}\right) (x).
\end{eqnarray*}%
And%
\begin{eqnarray*}
\left( \mu _{B_{1}}^{-}\circ \mu _{B_{2}}^{-}\right) (x) &=&{\bigvee }%
_{x=yz}\left\{ \mu _{B_{1}}^{-}\left( y\right) \vee \mu _{B_{2}}^{-}\left(
z\right) \right\}  \\
&\geq &{\bigvee }_{x=yz}\left\{ \mu _{B_{1}}^{-}\left( yz\right) \vee \mu
_{B_{2}}^{-}\left( yz\right) \right\}  \\
&=&\left\{ \mu _{B_{1}}^{-}\left( x\right) \vee \mu _{B_{2}}^{-}\left(
x\right) \right\}  \\
&=&\left( \mu _{B_{1}}^{-}\cup \mu _{B_{2}}^{-}\right) (x).
\end{eqnarray*}%
Thus we get $\mu _{B_{1}}^{+}\circ \mu _{B_{2}}^{+}\subseteq \mu
_{B_{1}}^{+}\cap \mu _{B_{2}}^{+}$ and $\mu _{B_{1}}^{-}\circ \mu
_{B_{2}}^{-}\supseteq \mu _{B_{1}}^{-}\cup \mu _{B_{2}}^{-}.$ Hence $%
B_{1}\circ B_{2}\subseteq B_{1}\cap B_{2}.$ $\ \ \Box $

\begin{proposition}
\label{P101}\textit{Let }$B_{1}=\left\langle \mu _{B_{1}}^{+},\mu
_{B_{1}}^{-}\right\rangle $\textit{\ and }$B_{2}=\left\langle \mu
_{B_{2}}^{+},\mu _{B_{2}}^{-}\right\rangle $\textit{\ be two
BVF-LA-subsemigroups of }$S$\textit{. Then }$B_{1}\cap B_{2}$\textit{\ is
also a BVF-LA-subsemigroup of }$S$\textit{.\ }
\end{proposition}

\textbf{Proof.\ }Let $B_{1}=\left\langle \mu _{B_{1}}^{+},\mu
_{B_{1}}^{-}\right\rangle $ and $B_{2}=\left\langle \mu _{B_{2}}^{+},\mu
_{B_{2}}^{-}\right\rangle $ be two \textit{BVF-LA-subsemigroups} of $S$. Let 
$x,y\in S$. Then%
\begin{eqnarray*}
\left( \mu _{B_{1}}^{+}\cap \mu _{B_{2}}^{+}\right) \left( xy\right)  &=&\mu
_{B_{1}}^{+}\left( xy\right) \wedge \mu _{B_{2}}^{+}\left( xy\right)  \\
&\geq &\left( \mu _{B_{1}}^{+}\left( x\right) \wedge \mu _{B_{1}}^{+}\left(
y\right) \right) \wedge \left( \mu _{B_{2}}^{+}\left( x\right) \wedge \mu
_{B_{2}}^{+}\left( y\right) \right)  \\
&=&\left( \mu _{B_{1}}^{+}\left( x\right) \wedge \mu _{B_{2}}^{+}\left(
x\right) \right) \wedge \left( \mu _{B_{1}}^{+}\left( y\right) \wedge \mu
_{B_{2}}^{+}\left( y\right) \right)  \\
&=&\left( \mu _{B_{1}}^{+}\cap \mu _{B_{2}}^{+}\right) \left( x\right)
\wedge \left( \mu _{B_{1}}^{+}\cap \mu _{B_{2}}^{+}\right) \left( y\right) .
\end{eqnarray*}%
And%
\begin{eqnarray*}
\left( \mu _{B_{1}}^{-}\cup \mu _{B_{2}}^{-}\right) \left( xy\right)  &=&\mu
_{B_{1}}^{-}\left( xy\right) \vee \mu _{B_{2}}^{-}\left( xy\right)  \\
&\leq &\left( \mu _{B_{1}}^{-}\left( x\right) \vee \mu _{B_{1}}^{-}\left(
y\right) \right) \vee \left( \mu _{B_{2}}^{-}\left( x\right) \vee \mu
_{B_{2}}^{-}\left( y\right) \right)  \\
&=&\left( \mu _{B_{1}}^{-}\left( x\right) \vee \mu _{B_{2}}^{-}\left(
x\right) \right) \vee \left( \mu _{B_{1}}^{-}\left( y\right) \vee \mu
_{B_{2}}^{-}\left( y\right) \right)  \\
&=&\left( \mu _{B_{1}}^{-}\cup \mu _{B_{2}}^{-}\right) \left( x\right) \vee
\left( \mu _{B_{1}}^{-}\cup \mu _{B_{2}}^{-}\right) \left( y\right) .
\end{eqnarray*}%
Thus $B_{1}\cap B_{2}$\textit{\ is also a bipolar-valued fuzzy
LA-subsemigroup of }$S$\textit{.} $\ \ \Box $

\begin{proposition}
\textit{Let }$B_{1}=\left\langle \mu _{B_{1}}^{+},\mu
_{B_{1}}^{-}\right\rangle $\textit{\ and }$B_{2}=\left\langle \mu
_{B_{2}}^{+},\mu _{B_{2}}^{-}\right\rangle $\textit{\ be two BVF-left (resp.
BVF-right, BVF-two-sided) ideal of }$S$\textit{. Then }$B_{1}\cap B_{2}$%
\textit{\ is also a BVF-left (resp. BVF-right, BVF-two-sided) ideal of }$S.$
\end{proposition}

\textbf{Proof.\ }The proof is similar to the proof of Proposition \ref{P101}%
. $\ \ \Box $

\begin{lemma}
\label{L2}In an LA-semigroup $S$ with left identity, for every BVF-left
ideal $B=\left\langle \mu _{B}^{+},\mu _{B}^{-}\right\rangle $ of $S$, we
have $\Gamma \circ B=B.$ Where $\Gamma =\left\langle \mathcal{S}_{\Gamma
}^{+}(x),\mathcal{S}_{\Gamma }^{-}(x)\right\rangle .$
\end{lemma}

\textbf{Proof. }Let $B=\left\langle \mu _{B}^{+},\mu _{B}^{-}\right\rangle $
be a BVF-left ideal of $S.$ It is sufficient to show that $\mathcal{S}%
_{\Gamma }^{+}\circ \mu _{B}^{+}\subseteq \mu _{B}^{+}$ and $\mathcal{S}%
_{\Gamma }^{-}\circ \mu _{B}^{-}\supseteq \mu _{B}^{-}$. Now $x=ex$, for all 
$x$ in $S$, as $e$ is left identity in $S$. So%
\begin{equation*}
(\mathcal{S}_{\Gamma }^{+}\circ \mu _{B}^{+})(x)={\bigvee }_{x=yz}\left\{ 
\mathcal{S}_{\Gamma }^{+}(y)\wedge \mu _{B}^{+}(z)\right\} \geq \mathcal{S}%
_{\Gamma }^{+}(e)\wedge \mu _{B}^{+}(x)=1\wedge \mu _{B}^{+}(x)=\mu
_{B}^{+}(x).
\end{equation*}%
And%
\begin{equation*}
(\mathcal{S}_{\Gamma }^{-}\circ \mu _{B}^{-})(x)={\bigwedge }_{x=yz}\left\{ 
\mathcal{S}_{\Gamma }^{-}(y)\vee \mu _{B}^{-}(z)\right\} \leq \mathcal{S}%
_{\Gamma }^{-}(e)\vee \mu _{B}^{-}(x)=-1\vee \mu _{B}^{-}(x)=\mu _{B}^{-}(x).
\end{equation*}%
Thus $\mathcal{S}_{\Gamma }^{+}\circ \mu _{B}^{+}\supseteq \mu _{B}^{+}$ and 
$\mathcal{S}_{\Gamma }^{-}\circ \mu _{B}^{-}\subseteq \mu _{B}^{-}.$ Hence $%
\Gamma \circ B=B.$ $\ \ \Box $

\textbf{Definition 3.4 }Let $S$ be an LA-semigroup and let $\emptyset \neq
A\subseteq S.$ Then bipolar-valued fuzzy characteristic function $\chi
_{A}=\left\langle \mu _{\chi _{A}}^{+},\mu _{\chi _{A}}^{-}\right\rangle $
of $A$ is defined as%
\begin{equation*}
\mu _{\chi _{A}}^{+}=\left\{ 
\begin{array}{c}
1\text{ \ \ if }x\in A \\ 
0\text{ \ \ if }x\notin A%
\end{array}%
\right. \text{ \ \ and \ }\mu _{\chi _{A}}^{-}=\left\{ 
\begin{array}{c}
-1\text{ \ \ if }x\in A \\ 
\text{ \ }0\text{ \ \ if }x\notin A.%
\end{array}%
\right.
\end{equation*}

\begin{theorem}
\label{T119}Let $A$ be a nonempty subset of an LA-semigroup $S$. Then $A$ is
an LA-subsemigroup of $S$ if and only if $\chi _{A}$ is a
BVF-LA-subsemigroup of $S$.
\end{theorem}

\textbf{Proof. }Let $A$ be an LA-subsemigroup of $S$. For any $x,y\in S,$ we
have the following cases:

Case $\left( 1\right) :$ If $x,y\in A$, then $xy\in A$. Since $A$ is an
LA-subsemigroup of $S$. Then $\mu _{\chi _{A}}^{+}\left( xy\right) =1,$ $\mu
_{\chi _{A}}^{+}\left( x\right) =1$ and $\mu _{\chi _{A}}^{+}\left( y\right)
=1$. Therefore 
\begin{equation*}
\mu _{\chi _{A}}^{+}\left( xy\right) =\mu _{\chi _{A}}^{+}\left( x\right)
\wedge \mu _{\chi _{A}}^{+}\left( y\right) .
\end{equation*}%
And $\mu _{\chi _{A}}^{-}\left( xy\right) =-1,$ $\mu _{\chi _{A}}^{-}\left(
x\right) =-1$ and $\mu _{\chi _{A}}^{-}\left( y\right) =-1$. Therefore 
\begin{equation*}
\mu _{\chi _{A}}^{-}\left( xy\right) =\mu _{\chi _{A}}^{-}\left( x\right)
\vee \mu _{\chi _{A}}^{-}\left( y\right) .
\end{equation*}%
Case $\left( 2\right) :$ If $x,y\notin A$, then $\mu _{\chi _{A}}^{+}\left(
x\right) =0$ and $\mu _{\chi _{A}}^{+}\left( y\right) =0$. So \ \ \ \ \ \ \ 
\begin{equation*}
\text{ \ \ \ \ \ \ \ }\mu _{\chi _{A}}^{+}\left( xy\right) \geq 0=\mu _{\chi
_{A}}^{+}\left( x\right) \wedge \mu _{\chi _{A}}^{+}\left( y\right) .
\end{equation*}%
And $\mu _{\chi _{A}}^{-}\left( x\right) =0$ and $\mu _{\chi _{A}}^{-}\left(
y\right) =0$. So 
\begin{equation*}
\text{ \ \ \ \ \ \ }\mu _{\chi _{A}}^{-}\left( xy\right) \leq 0=\mu _{\chi
_{A}}^{-}\left( x\right) \vee \mu _{\chi _{A}}^{-}\left( y\right) .
\end{equation*}

Case $\left( 3\right) :$ If $x\in A$ or $y\in A$. If $x\in A$ and $y\notin A$%
, then $\mu _{\chi _{A}}^{+}\left( x\right) =1$ and $\mu _{\chi
_{A}}^{+}\left( y\right) =0$. So%
\begin{equation*}
\text{ \ \ \ \ \ \ \ \ }\mu _{\chi _{A}}^{+}\left( xy\right) \geq 0=\mu
_{\chi _{A}}^{+}\left( x\right) \wedge \mu _{\chi _{A}}^{+}\left( y\right) .
\end{equation*}%
Now if $x\notin A$ and $y\in A$, then $\mu _{\chi _{A}}^{+}\left( x\right) =0
$ and $\mu _{\chi _{A}}^{+}\left( y\right) =1$. So%
\begin{equation*}
\text{ \ \ \ \ \ \ \ \ \ }\mu _{\chi _{A}}^{+}\left( xy\right) \geq 0=\mu
_{\chi _{A}}^{+}\left( x\right) \wedge \mu _{\chi _{A}}^{+}\left( y\right) .
\end{equation*}%
And if $x\in A$ or $y\in A$. If $x\in A$ and $y\notin A$, then $\mu _{\chi
_{A}}^{-}\left( x\right) =-1$ and $\mu _{\chi _{A}}^{-}\left( y\right) =0$.
So%
\begin{equation*}
\text{ \ \ \ \ \ \ \ }\mu _{\chi _{A}}^{-}\left( xy\right) \leq 0=\mu _{\chi
_{A}}^{-}\left( x\right) \vee \mu _{\chi _{A}}^{-}\left( y\right) .
\end{equation*}%
Now if If $x\notin A$ and $y\in A$, then $\mu _{\chi _{A}}^{-}\left(
x\right) =0$ and $\mu _{\chi _{A}}^{-}\left( y\right) =-1$. So%
\begin{equation*}
\text{ \ \ \ \ \ \ \ \ }\mu _{\chi _{A}}^{-}\left( xy\right) \leq 0=\mu
_{\chi _{A}}^{-}\left( x\right) \vee \mu _{\chi _{A}}^{-}\left( y\right) .
\end{equation*}%
Hence $\chi _{A}=\left\langle \mu _{\chi _{A}}^{+},\mu _{\chi
_{A}}^{-}\right\rangle $ is a BVF-LA-subsemigroup of $S$.

Conversely, suppose $\chi _{A}=\left\langle \mu _{\chi _{A}}^{+},\mu _{\chi
_{A}}^{-}\right\rangle $ is a BVF-LA-subsemigroup of $S$ and let $x,y\in A$.
Then we have 
\begin{eqnarray*}
\text{ \ \ }\mu _{\chi _{A}}^{+}\left( xy\right) &\geq &\mu _{\chi
_{A}}^{+}\left( x\right) \wedge \mu _{\chi _{A}}^{+}\left( y\right) =1\wedge
1=1 \\
\mu _{\chi _{A}}^{+}\left( xy\right) &\geq &1\text{ but }\mu _{\chi
_{A}}^{+}\left( xy\right) \leq 1 \\
\mu _{\chi _{A}}^{+}\left( xy\right) &=&1
\end{eqnarray*}%
And 
\begin{eqnarray*}
\text{ \ \ }\mu _{\chi _{A}}^{-}\left( xy\right) &\leq &\mu _{\chi
_{A}}^{-}\left( x\right) \vee \mu _{\chi _{A}}^{-}\left( y\right) =-1\vee
-1=-1 \\
\mu _{\chi _{A}}^{-}\left( xy\right) &\leq &-1\text{ but }\mu _{\chi
_{A}}^{-}\left( xy\right) \geq -1 \\
\mu _{\chi _{A}}^{-}\left( xy\right) &=&-1
\end{eqnarray*}%
Hence $xy\in A$. Therefore $A$ is an LA-subsemigroup of $S.$ $\ \ \Box $

\begin{theorem}
Let $A$ be a nonempty subset of an LA-semigroup $S$. Then $A$ is a left
(resp. right) ideal of $S$ if and only if $\chi _{A}$ is a BVF-left (resp.
BVF-right) ideal of $S$.
\end{theorem}

\textbf{Proof. }The proof of this theorem is similar to Theorem \ref{T119}. $%
\ \ \Box $

\textbf{Definition 3.5 }A BVF-subset $B=\left\langle \mu _{B}^{+},\mu
_{B}^{-}\right\rangle $ of an LA-semigroup $S$ is called a\textbf{\ }%
BVF-generalized bi-ideal\textbf{\ }of $S$ if%
\begin{equation*}
\mu _{B}^{+}\left( (xy)z\right) \geq \mu _{B}^{+}(x)\wedge \mu _{B}^{+}(y)%
\text{ \ and \ }\mu _{B}^{-}\left( (xy)z\right) \leq \mu _{B}^{-}(x)\vee \mu
_{B}^{-}(y),\text{ for all }x,y,z\in S.
\end{equation*}

\textbf{Definition 3.6 }A BVF-LA-subsemigroup $B=\left\langle \mu
_{B}^{+},\mu _{B}^{-}\right\rangle $ of an LA-semigroup $S$ is called a%
\textbf{\ }BVF-bi-ideal\textbf{\ }of $S$ if%
\begin{equation*}
\mu _{B}^{+}\left( (xy)z\right) \geq \mu _{B}^{+}(x)\wedge \mu _{B}^{+}(y)%
\text{ \ and \ }\mu _{B}^{-}\left( (xy)z\right) \leq \mu _{B}^{-}(x)\vee \mu
_{B}^{-}(y),\text{ for all }x,y,z\in S.
\end{equation*}

\begin{lemma}
\label{L57}A BVF-subset $B=\left\langle \mu _{B}^{+},\mu
_{B}^{-}\right\rangle $ of an LA-semigroup $S$ is a BVF-generalized bi-ideal
of $S$ if and only if $\left( \mu _{B}^{+}\circ \mathcal{S}_{\Gamma
}^{+}\right) \circ \mu _{B}^{+}\subseteq \mu _{B}^{+}$ and $\left( \mu
_{B}^{-}\circ \mathcal{S}_{\Gamma }^{-}\right) \circ \mu _{B}^{-}\supseteq
\mu _{B}^{-}.$
\end{lemma}

\textbf{Proof. }Let $B=\left\langle \mu _{B}^{+},\mu _{B}^{-}\right\rangle $
be a BVF-generalized bi-ideal of an LA-semigroup $S$ and $x\in S.$ If $%
\left( \left( \mu _{B}^{+}\circ \mathcal{S}_{\Gamma }^{+}\right) \circ \mu
_{B}^{+}\right) \left( x\right) =0$ and $\left( \left( \mu _{B}^{-}\circ 
\mathcal{S}_{\Gamma }^{-}\right) \circ \mu _{B}^{-}\right) (x)=0,$ then%
\begin{equation*}
\left( \left( \mu _{B}^{+}\circ \mathcal{S}_{\Gamma }^{+}\right) \circ \mu
_{B}^{+}\right) \left( x\right) =0\leq \mu _{B}^{+}\left( x\right) \text{ \
and \ }\left( \left( \mu _{B}^{-}\circ \mathcal{S}_{\Gamma }^{-}\right)
\circ \mu _{B}^{-}\right) (x)=0\geq \mu _{B}^{-}(x).
\end{equation*}%
Otherwise%
\begin{eqnarray*}
\left( \left( \mu _{B}^{+}\circ \mathcal{S}_{\Gamma }^{+}\right) \circ \mu
_{B}^{+}\right) \left( x\right) &=&{\bigvee }_{x=ab}\left\{ \left( \mu
_{B}^{+}\circ \mathcal{S}_{\Gamma }^{+}\right) \left( a\right) \wedge \mu
_{B}^{+}\left( b\right) \right\} \\
&=&{\bigvee }_{x=ab}\left\{ {\bigvee }_{a=mn}\left\{ \mu _{B}^{+}\left(
m\right) \wedge \mathcal{S}_{\Gamma }^{+}\left( n\right) \right\} \wedge \mu
_{B}^{+}\left( b\right) \right\} \\
&=&{\bigvee }_{x=ab}{\bigvee }_{a=mn}\left\{ \left( \mu _{B}^{+}\left(
m\right) \wedge 1\right) \wedge \mu _{B}^{+}\left( b\right) \right\} \\
&=&{\bigvee }_{x=ab}{\bigvee }_{a=mn}\left\{ \mu _{B}^{+}\left( m\right)
\wedge \mu _{B}^{+}\left( b\right) \right\} \\
&\leq &\mu _{B}^{+}\left( x\right) .
\end{eqnarray*}%
And%
\begin{eqnarray*}
\left( \left( \mu _{B}^{-}\circ \mathcal{S}_{\Gamma }^{-}\right) \circ \mu
_{B}^{-}\right) \left( x\right) &=&{\bigwedge }_{x=ab}\left\{ \left( \mu
_{B}^{-}\circ \mathcal{S}_{\Gamma }^{-}\right) \left( a\right) \vee \mu
_{B}^{-}\left( b\right) \right\} \\
&=&{\bigwedge }_{x=ab}\left\{ {\bigwedge }_{a=mn}\left\{ \mu _{B}^{-}\left(
m\right) \vee \mathcal{S}_{\Gamma }^{-}\left( n\right) \right\} \vee \mu
_{B}^{-}\left( b\right) \right\} \\
&=&{\bigwedge }_{x=ab}{\bigwedge }_{a=mn}\left\{ \left( \mu _{B}^{-}\left(
m\right) \vee -1\right) \vee \mu _{B}^{-}\left( b\right) \right\} \\
&=&{\bigwedge }_{x=ab}{\bigwedge }_{a=mn}\left\{ \mu _{B}^{-}\left( m\right)
\vee \mu _{B}^{-}\left( b\right) \right\} \\
&\geq &\mu _{B}^{-}\left( x\right) .
\end{eqnarray*}%
Thus $\left( \mu _{B}^{+}\circ \mathcal{S}_{\Gamma }^{+}\right) \circ \mu
_{B}^{+}\subseteq \mu _{B}^{+}$ and $\left( \mu _{B}^{-}\circ \mathcal{S}%
_{\Gamma }^{-}\right) \circ \mu _{B}^{-}\supseteq \mu _{B}^{-}.$\newline
Conversely, assume that $\left( \mu _{B}^{+}\circ \mathcal{S}_{\Gamma
}^{+}\right) \circ \mu _{B}^{+}\subseteq \mu _{B}^{+}$ and $\left( \mu
_{B}^{-}\circ \mathcal{S}_{\Gamma }^{-}\right) \circ \mu _{B}^{-}\supseteq
\mu _{B}^{-}.$ Let $x,y,z\in S,$ then%
\begin{eqnarray*}
\mu _{B}^{+}\left( (xy)z\right) &\geq &\left( \left( \mu _{B}^{+}\circ 
\mathcal{S}_{\Gamma }^{+}\right) \circ \mu _{B}^{+}\right) \left(
(xy)z\right) ={\bigvee }_{(xy)z=cd}\left\{ \left( \mu _{B}^{+}\circ \mathcal{%
S}_{\Gamma }^{+}\right) \left( c\right) \wedge \mu _{B}^{+}\left( d\right)
\right\} \\
&\geq &\left( \mu _{B}^{+}\circ \mathcal{S}_{\Gamma }^{+}\right) \left(
xy\right) \wedge \mu _{B}^{+}\left( z\right) =\left\{ {\bigvee }%
_{xy=pq}\left\{ \mu _{B}^{+}\left( p\right) \wedge \mathcal{S}_{\Gamma
}^{+}\left( q\right) \right\} \right\} \wedge \mu _{B}^{+}\left( z\right) \\
&\geq &\left\{ \mu _{B}^{+}\left( x\right) \wedge \mathcal{S}_{\Gamma
}^{+}\left( y\right) \right\} \wedge \mu _{B}^{+}\left( z\right) =\left\{
\mu _{B}^{+}\left( x\right) \wedge 1\right\} \wedge \mu _{B}^{+}\left(
z\right) \\
&=&\mu _{B}^{+}\left( x\right) \wedge \mu _{B}^{+}\left( z\right) .
\end{eqnarray*}%
And%
\begin{eqnarray*}
\mu _{B}^{-}\left( (xy)z\right) &\leq &\left( \left( \mu _{B}^{-}\circ 
\mathcal{S}_{\Gamma }^{-}\right) \circ \mu _{B}^{-}\right) \left(
(xy)z\right) ={\bigwedge }_{(xy)z=cd}\left\{ \left( \mu _{B}^{-}\circ 
\mathcal{S}_{\Gamma }^{-}\right) \left( c\right) \vee \mu _{B}^{-}\left(
d\right) \right\} \\
&\leq &\left( \mu _{B}^{-}\circ \mathcal{S}_{\Gamma }^{-}\right) \left(
xy\right) \vee \mu _{B}^{-}\left( z\right) =\left\{ {\bigwedge }%
_{xy=pq}\left\{ \mu _{B}^{-}\left( p\right) \vee \mathcal{S}_{\Gamma
}^{-}\left( q\right) \right\} \right\} \vee \mu _{B}^{-}\left( z\right) \\
&\leq &\left\{ \mu _{B}^{-}\left( x\right) \vee \mathcal{S}_{\Gamma
}^{-}\left( y\right) \right\} \vee \mu _{B}^{-}\left( z\right) =\left\{ \mu
_{B}^{-}\left( x\right) \vee -1\right\} \vee \mu _{B}^{-}\left( z\right) \\
&=&\mu _{B}^{-}\left( x\right) \vee \mu _{B}^{-}\left( z\right) .
\end{eqnarray*}%
Thus $\mu _{B}^{+}\left( (xy)z\right) \geq \mu _{B}^{+}\left( x\right)
\wedge \mu _{B}^{+}\left( z\right) $ and $\mu _{B}^{-}\left( (xy)z\right)
\leq \mu _{B}^{-}\left( x\right) \vee \mu _{B}^{-}\left( z\right) ,$ which
implies that $B=\left\langle \mu _{B}^{+},\mu _{B}^{-}\right\rangle $ is a
BVF-generalized bi-ideal of $S.$ $\ \ \Box $

\begin{lemma}
\label{L30}Let $B=\left\langle \mu _{B}^{+},\mu _{B}^{-}\right\rangle $ be a
BVF-subset of an LA-semigroup $S$ then $B=\left\langle \mu _{B}^{+},\mu
_{B}^{-}\right\rangle $ is a BVF-bi-ideal of $S$\ if and only if $\mu
_{B}^{+}\circ \mu _{B}^{+}\subseteq \mu _{B}^{+},$ $\mu _{B}^{-}\circ \mu
_{B}^{-}\supseteq \mu _{B}^{-},$ $\left( \mu _{B}^{+}\circ \mathcal{S}%
_{\Gamma }^{+}\right) \circ \mu _{B}^{+}\subseteq \mu _{B}^{+}$ and $\left(
\mu _{B}^{-}\circ \mathcal{S}_{\Gamma }^{-}\right) \circ \mu
_{B}^{-}\supseteq \mu _{B}^{-}.$
\end{lemma}

\textbf{Proof. }Follows from Lemma \ref{L56}(1) and Lemma \ref{L57}. $\ \
\Box $

\textbf{Definition 3.7 }A BVF-subset $B=\left\langle \mu _{B}^{+},\mu
_{B}^{-}\right\rangle $ of an LA-semigroup $S$ is called a BVF-interior
ideal of $S$ if%
\begin{equation*}
\mu _{B}^{+}\left( (xy)z\right) \geq \mu _{B}^{+}\left( y\right) \text{ \
and \ }\mu _{B}^{-}\left( (xy)z\right) \leq \mu _{B}^{-}\left( y\right) ,%
\text{ for all }x,y,z\in S.
\end{equation*}

\begin{lemma}
\label{L31}Let $B=\left\langle \mu _{B}^{+},\mu _{B}^{-}\right\rangle $ be a
BVF-subset of an LA-semigroup $S$ then $B=\left\langle \mu _{B}^{+},\mu
_{B}^{-}\right\rangle $ is a BVF-interior ideal of $S$\ if and only if $%
\left( \mathcal{S}_{\Gamma }^{+}\circ \mu _{B}^{+}\right) \circ \mathcal{S}%
_{\Gamma }^{+}\subseteq \mu _{B}^{+}$ and $\left( \mathcal{S}_{\Gamma
}^{-}\circ \mu _{B}^{-}\right) \circ \mathcal{S}_{\Gamma }^{-}\supseteq \mu
_{B}^{-}.$
\end{lemma}

\textbf{Proof. }The proof of this lemma is similar to the proof of Lemma \ref%
{L57}$.$ $\ \ \Box $

\begin{remark}
Every BVF-ideal is a BVF-interior ideal of an LA-semigroup $S,$ but the
converse is not true.
\end{remark}

\textbf{Example 3.2 }Let $S=\{a,b,c,d\}$, the binary operation "$\cdot $" on 
$S$ be defined as follows:%
\begin{equation*}
\begin{tabular}{l|llll}
$\cdot $ & $a$ & $b$ & $c$ & $d$ \\ \hline
$a$ & $c$ & $c$ & $c$ & $d$ \\ 
$b$ & $d$ & $d$ & $c$ & $c$ \\ 
$c$ & $d$ & $d$ & $d$ & $d$ \\ 
$d$ & $d$ & $d$ & $d$ & $d$%
\end{tabular}%
\end{equation*}%
Clearly, $S$ is an LA-semigroup. But $S$ is not a semigroup because $%
c=a\cdot (a\cdot b)\neq (a\cdot a)\cdot b=d.$ Now we define BVF-subset as%
\begin{equation*}
B=\left\langle \mu _{B}^{+},\mu _{B}^{-}\right\rangle =\left\langle \left( 
\frac{a}{0.5},\frac{b}{0.3},\frac{c}{0.1},\frac{d}{0.8}\right) ,\text{ }%
\left( \frac{a}{-0.7},\frac{b}{-0.4},\frac{c}{-0.2},\frac{d}{-0.9}\right)
\right\rangle .
\end{equation*}%
It can be verified that $B=\left\langle \mu _{B}^{+},\mu
_{B}^{-}\right\rangle $ is a BVF-interior ideal of $S.$ But, since%
\begin{equation*}
\mu _{B}^{+}\left( b\cdot c\right) =\mu _{B}^{+}\left( c\right) =0.1<0.3=\mu
_{B}^{+}\left( b\right) .
\end{equation*}%
And%
\begin{equation*}
\mu _{B}^{-}\left( b\cdot c\right) =\mu _{B}^{-}\left( c\right)
=-0.2>-0.4=\mu _{B}^{-}\left( b\right) .
\end{equation*}%
Thus $B=\left\langle \mu _{B}^{+},\mu _{B}^{-}\right\rangle $ is not a
BVF-right ideal of $S,$ that is, $B=\left\langle \mu _{B}^{+},\mu
_{B}^{-}\right\rangle $ is not a BVF-two-sided ideal of $S.$

\begin{proposition}
Every BVF-subset $B=\left\langle \mu _{B}^{+},\mu _{B}^{-}\right\rangle $ of
an LA-semigroup $S$ with left identity is a BVF-right ideal if and only if
it is a BVF-interior ideal.
\end{proposition}

\textbf{Proof. }Let every BVF-subset $B=\left\langle \mu _{B}^{+},\mu
_{B}^{-}\right\rangle $ of $S$ is a BVF-right ideal. For $x$, $a$ and $y$ of 
$S$, consider%
\begin{equation*}
\mu _{B}^{+}((xa)y)\geq \mu _{B}^{+}(xa)=\mu _{B}^{+}((ex)a)=\mu
_{B}^{+}((ax)e)\geq \mu _{B}^{+}(ax)\geq \mu _{B}^{+}(a).
\end{equation*}%
And%
\begin{equation*}
\mu _{B}^{-}((xa)y)\leq \mu _{B}^{-}(xa)=\mu _{B}^{-}((ex)a)=\mu
_{B}^{-}((ax)e)\leq \mu _{B}^{-}(ax)\leq \mu _{B}^{-}(a).
\end{equation*}%
Which implies that $B=\left\langle \mu _{B}^{+},\mu _{B}^{-}\right\rangle $
is a BVF-interior ideal. Conversely, for any $x$ and $y$ in $S$ we have,%
\begin{equation*}
\mu _{B}^{+}(xy)=\mu _{B}^{+}((ex)y)\geq \mu _{B}^{+}(x)\text{ \ and \ }\mu
_{B}^{-}(xy)=\mu _{B}^{-}((ex)y)\leq \mu _{B}^{-}(x).
\end{equation*}%
Hence required. $\ \ \Box $

\begin{theorem}
Let $B=\left\langle \mu _{B}^{+},\mu _{B}^{-}\right\rangle $ be a BVF-left
ideal of an LA-semigroup $S$ with left identity, then $B=\left\langle \mu
_{B}^{+},\mu _{B}^{-}\right\rangle $ being BVF-interior ideal is a
BVF-bi-ideal of $S$.
\end{theorem}

\textbf{Proof. }Since $B=\left\langle \mu _{B}^{+},\mu _{B}^{-}\right\rangle 
$ is an BVF-left ideal in $S$, so $\mu _{B}^{+}(xy)\geq \mu _{B}^{+}(y)$ and 
$\mu _{B}^{-}(xy)\leq \mu _{B}^{-}(y)$ for all $x$ and $y$ in $S$. As $e$ is
left identity in $S$. So,%
\begin{equation*}
\mu _{B}^{+}(xy)=\mu _{B}^{+}((ex)y)\geq \mu _{B}^{+}(x)\text{ \ and \ }\mu
_{B}^{-}(xy)=\mu _{B}^{-}((ex)y)\leq \mu _{B}^{-}(x),
\end{equation*}
which implies that $\mu _{B}^{+}(xy)\geq \mu _{B}^{+}(x)\wedge \mu
_{B}^{+}(y)$ and $\mu _{B}^{-}(xy)\leq \mu _{B}^{-}(x)\vee \mu _{B}^{-}(y)$
for all $x$ and $y$ in $S$. Thus $B=\left\langle \mu _{B}^{+},\mu
_{B}^{-}\right\rangle $ is an BVF-LA-subsemigroup of $S$. For any $x,y$ and $%
z$ in $S$, we get%
\begin{equation*}
\mu _{B}^{+}((xy)z)=\mu _{B}^{+}((x(ey))z)=\mu _{B}^{+}((e(xy))z)\geq \mu
_{B}^{+}(xy)=\mu _{B}^{+}((ex)y)\geq \mu _{B}^{+}(x).
\end{equation*}%
And%
\begin{equation*}
\mu _{B}^{-}((xy)z)=\mu _{B}^{-}((x(ey))z)=\mu _{B}^{-}((e(xy))z)\leq \mu
_{B}^{-}(xy)=\mu _{B}^{-}((ex)y)\leq \mu _{B}^{-}(x).
\end{equation*}%
Also%
\begin{equation*}
\mu _{B}^{+}((xy)z)=\mu _{B}^{+}((zy)x)=\mu _{B}^{+}((z(ey))x)=\mu
_{B}^{+}((e(zy))x)\geq \mu _{B}^{+}(zy)=\mu _{B}^{+}((ez)y)\geq \mu
_{B}^{+}(z).
\end{equation*}%
And%
\begin{equation*}
\mu _{B}^{-}((xy)z)=\mu _{B}^{-}((zy)x)=\mu _{B}^{-}((z(ey))x)=\mu
_{B}^{-}((e(zy))x)\leq \mu _{B}^{-}(zy)=\mu _{B}^{-}((ez)y)\leq \mu
_{B}^{-}(z).
\end{equation*}%
Hence $\mu _{B}^{+}((xy)z)\geq \mu _{B}^{+}(x)\wedge \mu _{B}^{+}(z)$ and $%
\mu _{B}^{-}((xy)z)\leq \mu _{B}^{-}(x)\vee \mu _{B}^{-}(z)$ for all $x,y$
and $z$ in $S$. $\ \ \Box $

\end{document}